Optimal control of second-order integral equations

S. A. Belbas.  SBELBAS@GMAIL.COM
Dept. of Mathematics, University of Alabama, Tuscaloosa, AL 35487-0350, USA.

Abstract. We analyze optimal control problems for multiple Fredholm and Volterra integral equations. These are nom-Pontryaginian optimal control problems, i.e. an extremum principle of Pontryagin's type does not hold. We obtain first-order necessary conditions for optimality, and second-order necessary and sufficient conditions. We illustrate with applications to first – and second – order Volterra bilinear control problems.

1. Double integral equations.

In this paper we introduce a new category of control stems, viz. systems governed by double integral equations. "Double" refers to multiple integrations, with multiplicity twice the dimension of the independent variable. It is relatively straightforward to extend the results and methods to higher degree of multiplicity. The crucial transition if from first-order integral equations (with integrals of the same multiplicity as the independent variable) to second-order integral systems.

The two basic types of integral equations of the second kind, Fredholm and Volterra, are well known. A (generally nonlinear) Fredholm integral equation of the second kind has the form

$$\varphi(x) = \varphi_0(x) + \int_G f(x, y, \varphi(y)) \, dy \tag{1.1}$$

and a Volterra integral equation of the second type has the general form

$$y(t) = y_0(t) + \int_0^t f(t, s, y(s)) \, ds \tag{1.2}$$

Here, we shall be interested in two types of double integral equations of the second kind. The two types (Fredholm and Volterra) of double integral equations are

$$\varphi(x) = \varphi_0(x) + \int_G f_1(x, y, \varphi(y)) \, dy + \tfrac{1}{2} \int_G \int_G f_2(x, y, z, \varphi(y), \varphi(z)) \, dy \, dz \, ;$$

$$y(t) = y_0(t) + \int_0^t f_1(t, s, y(s)) \, ds + \tfrac{1}{2} \int_0^t \int_0^t f_2(t, s, \sigma, y(s), y(\sigma)) \, d\sigma \, ds \tag{1.3}$$

These double integral equations are particular cases of the general multiple integral equations

$$\varphi(x) = \varphi_0(x) + \sum_{k=1}^{\infty} \frac{1}{k!} \int_G \int_G \cdots \int_G f_k(x, y_1, y_2, \ldots, y_k, \varphi(y_1), \varphi(y_2), \ldots, \varphi(y_k)) \, dy_1 \, dy_2 \cdots dy_k \, ;$$

$$y(t) = y_0(t) + \sum_{k=1}^{\infty} \frac{1}{k!} \int_0^t \int_0^t \cdots \int_0^t f_k(t, s_1, s_2, \ldots, s_k, y(s_1), y(s_2), \ldots, y(s_k)) \, ds_1 \, ds_2 \cdots ds_k \tag{1.4}$$

In the terminology we are using here, *multiple integral equation* refers to multiple integrations of multiplicity higher than the dimension of the independent variable. Particular cases of multiple integral equation were introduced by E. Schmidt [S] for Fredholm integral equations, and by Tonelli [T] for Volterra integral equations. For general causal systems, potentially including multiple and infinite-order Volterra integral equations, [C] is a crucial reference. General multiple and infinite-order Volterra equations have been analyzed in [BB].

The control systems associated with the two types of multiple integral equations have the form

$$\varphi(x) = \varphi_0(x) + \sum_{k=1}^{\infty} \frac{1}{k!} \int_G \int_G \cdots \int_G f_k(x, y_1, y_2, ..., y_k, \varphi(y_1), \varphi(y_2), ..., \varphi(y_k),$$

$$u(y_1), u(y_2), ..., u(y_k)) \, dy_1 \, dy_2 \cdots dy_k \, ;$$

$$y(t) = y_0(t) + \sum_{k=1}^{\infty} \frac{1}{k!} \int_0^t \int_0^t \cdots \int_0^t f_k(t, s_1, s_2, ..., s_k, y(s_1), y(s_2), ..., y(s_k),$$

$$u(s_1), u(s_2), ..., u(s_k)) \, ds_1 \, ds_2 \cdots ds_k$$

(1.5)

where $u$ is a control function.

<u>Example 1.1.</u> We consider an integro-differential equation

$$\frac{dx}{dt} = f(t, x(t)) + \int_0^t g(t, \sigma, x(t), x(\sigma)) \, d\sigma \, ; \quad x(0) = x_0 \, .$$

(1.6)

We rewrite this equation in integral form:

$$x(t) = x_0 + \int_0^t f(s, x(s)) \, ds + \int_0^t \int_0^s g(s, \sigma, x(s), x(\sigma)) \, d\sigma \, ds \, .$$

(1.7)

We define the function

$$\tilde{g}(s, \sigma, x(s), x(\sigma)) := g(s, \sigma, x(s), x(\sigma)) + g(\sigma, s, x(\sigma), x(s)) \, .$$

(1.8)

Because the original function $g(t, \sigma, ...)$ is utilized (and therefore needs to be defined) only for $t \geq \sigma$, we have

$$\int_0^t \int_0^s g(s, \sigma, ...) \, d\sigma \, ds = \frac{1}{2} \left[ \int_0^t \int_0^s g(s, \sigma, ...) \, d\sigma \, ds + \int_0^t \int_0^\sigma g(\sigma, s, ...) \, ds \, d\sigma \right] =$$

$$= \frac{1}{2} \left[ \int_0^t \int_0^s g(s, \sigma, ...) \, d\sigma \, ds + \int_0^t \int_s^t g(\sigma, s ...) \, d\sigma \, ds \right] =$$

$$= \frac{1}{2} \int_0^t \int_0^t \tilde{g}(s, \sigma, ...) \, d\sigma \, ds$$

(1.9)

and consequently the integral equation becomes

$$x(t) = x_0 + \int_0^t f(s, x(s)) \, ds + \frac{1}{2} \int_0^t \int_0^t \tilde{g}(s, \sigma, x(s), x(\sigma)) \, d\sigma \, ds$$

(1.10)

The last integral equation can be further extended by introducing another ingredient of memory, making the functions appearing on the right-hand side also dependent on $t$, thus obtaining the integral equation

$$x(t) = x_0(t) + \int_0^t f(t,s,x(s))\,ds + \frac{1}{2}\int_0^t \int_0^t \tilde{g}(t,s,\sigma,x(s),x(\sigma))\,d\sigma\,ds .$$

(1.11)

<u>Example 1.2.</u> We consider a Voltera-Lotka integral system with controls corresponding to external controlled influx or harvesting of different components (for example, different biological species). We include distributed delays in the actions of the controls and in the interactions among different components. The system has the form

$$x(t) = x_0(t) + \int_0^t \left[A(t,s)x(s) + b^{jk}(t,s)x_j(s)x_k(s)\right]ds +$$

$$+ \frac{1}{2}\int_0^t \int_0^t \left[c^j(t,s,\sigma)x_j(s)u_j(\sigma) + d^{jk}(t,s,\sigma)x_j(s)x_k(\sigma)u_j(\sigma)u_k(\sigma)\right]d\sigma\,ds$$

(1.12)

In the above equation, the state $x$ takes values in $n$ – dimensional Euclidean space, and the control $u$ takes values in $m$ – dimensional Euclidean space. Summation with respect to repeated indices is used in the above system.

The double integral terms arise from the integral formulation of a Lotka integrodifferential model with distributed delay and instantaneous and delayed action of the control. This model is a natural extension of the basic integro-differential model of Lotka [L].

The system (1.12) represents an enhanced Volterra-Lotka model of interacting biogical populations $x_i(t)$, with $n$ biological species, and with harvesting controls representing rates of harvesting as percentages of the populations, and with distributed delays in the births and deaths/ inter-population interactions (predation), and the delayed effects of harvesting on birth rates. ///

We end this section with stating (for later use) the obvious fact that all results, in this paper, about second-order integral equations, can be specialized to first-order integral equations by the simple action of omitting some of the terms in the results proved in the present paper.

2.   <u>Multidimensional matrices.</u>

In this section, we gather certain definitions and introduce some notation, that will be useful in our treatment of second-order controlled integral equations.

We will be using <u>multidimensional</u> matrices, the term "multidimensional" referring to the dimensionality of the indices used for labelling the elements of the matrix. An ordinary matrix has two-dimensional indices, $A = \left[ A_{ij} : 1 \leq i \leq m, \, 1 \leq j \leq n \right]$; this matrix $A$ itself has dimensions $m$ by $n$; $m$ and $n$ represent the <u>range</u> of the indices $i$ and $j$; the <u>dimension</u> of $(i, j)$ is 2. When we want to differentiate between the two concepts of dimension, the dimension of a matrix and the dimension of the indices, we shall use the terms, respectively, <u>matricial dimension</u> and <u>indicial dimension</u>. The term <u>multidimensional</u> will refer to the indicial dimension of a matrix; with respect to matricial dimension, no special appellation is needed. A <u>multidimensional</u> matrix has indices of arbitrary dimension; a matrix of *indicial dimension $N$* has the form

$$A = \left[ A_{i_1 i_2 \cdots i_N} : \, 1 \leq i_k \leq m_k, \, 1 \leq k \leq N \right].$$

We use the term tridimensional matrix to denote an array indexed by 3 indices, for example $A_{ijk}$ or $A_k^{ij}$, and so on. Such arrays go also by the names "spatial matrices" or "matricial vectors", the latter terminology stemming from interpreting such a matrix as a vector whose components are matrices, e.g. $\left[ A_{ijk} \right]_{\substack{1 \leq i \leq \ell \\ 1 \leq j \leq m \\ 1 \leq k \leq n}}$ can be interpreted as an $n$ dimensional vector whose $k$ – th component is the $\ell \times m$ matrix $A_k := [A_{ijk}]_{\substack{1 \leq i \leq \ell \\ 1 \leq j \leq m}}$. (This interpretation is tacitly used, without attaching any specific appellation to it, in references on control of bilinear ordinary differential equations, e.g. [PY, D, M1, M2].)

References [SO1] and [SO2] represent a modern account of multidimensional matrices (in the sense of indicial dimension); the origins of multidimensional matrices go back to Hamilton (the same Hamilton of the Hamilton-Cayley theorem) and others; we omit detailed references of mainly historical relevance.

Some operations with tridimensional matrices, like addition and multiplication, are rather straightforward. There are several possible multiplications, for example, if

$A = \left[ A_{ijk} \right]$ and $B = \left[ B^{i'j'k'} \right]$, then, using the standard convention (of tensor algebra) of summation over an index that appears both as a subscript and as a superscript, in the same (matricial) monomial, we can define several products, for example

$\left[ A_{ijk} B^{ij'k'} \right], \left[ A_{ijk} B^{ijk'} \right]$ (provided the dimensions are compatible), etc. In each specific situation that involves products of arrays, either we shall specify the relevant product or we shall use explicit notation with subscripts and superscripts.

One aspect of tridimensional and multidimensional matrices is the concept of transposition (and taking adjoints, for matrices with complex entries). There are, of course, several different kinds of transformations that are akin to the ordinary transposition of ordinary (two-index) matrices. By a slight abuse of language, in order to avoid introducing too much new terminology, we will still call those transformations transpositions. For a tridimensional matrix $\left[ A_{ijk} \right]$, we take a rearrangement $\sigma$, of the symbols $i, j, k,$ and, for each such rearrangement, we define a transposition $T(\sigma)$ by

$$\left[ \left( A^{T(\sigma)} \right)_{ijk} \right] := \left[ A_{\sigma(i)\,\sigma(j)\,\sigma(k)} \right]$$. Analogous definitions can be given for multidimensional matrices involving subscripts and superscripts.

In general, for matrices with $n$ – dimensional indices,

$$\vec{i} = (i_1, i_2, ..., i_n)$$

if we also allow upper – and lower – indices (as we need to do for applications in control theory), there are $(n+1)!$ different matrices that can be formed, namely $n!$ ways of arranging the $n$ components of the multi-index $\vec{i}$ and, for each such arrangement, $(n+1)$ ways of splitting the arrangement into two ordered sets of upper – indices and lower – indices. Thus for every multi - indexed matrix, with distinction between upper – and lower – indices, there are $(n+1)!-1$ other matrices of the same indicial dimension that can be formed, thus $(n+1)!-1$ different transpositions (in the sense of matricial transposition).

A particular case of transpositions arises in the case of tridimensional matrices of the types $\left[ A_k^{ij} \right]$ and $\left[ A_{ij}^k \right]$. We define these for matrices of type $\left[ A_k^{ij} \right]$, and the definitions for matrices of the type $\left[ A_{ij}^k \right]$ follow by analogy. We place the symbols $i, j, k,$ in a triangular arrangement

$$\begin{matrix} i & & j \\ & k & \end{matrix}$$

, and we consider the rearrangements resulting from each of two rotations, viz. clockwise, leading to

$$\begin{matrix} k & & i \\ & j & \end{matrix}$$,

and counterclockwise, leading to

$j \qquad k$

$\qquad \cdot$

$i$

We denote these two rearrangements by $\sigma^-$ (for the clockwise rotation), and $\sigma^+$ (for the counterclockwise rotation). The two transpositions associated with these two rearrangements will be denoted by $T+$, $T-$, respectively, thus, in explicit notation,

$$\left[ \left( A^{T+} \right)^{ij}_k \right] := \left[ A^{ki}_j \right] , \quad \left[ \left( A^{T-} \right)^{ij}_k \right] := \left[ A^{jk}_i \right] .$$

A few more transpositions can be defined : the transpositions $T \updownarrow$, $T \leftrightarrow$, $T \searrow$, $T \nearrow$. These are defined by

$$\left( A^{T\updownarrow} \right)^{ij}_k = A^k_{ij} , \quad \left( A^{T\updownarrow} \right)^k_{ij} = A^{ij}_k ;$$

$$\left( A^{T\leftrightarrow} \right)^{ij}_k = A^{ji}_k , \quad \left( A^{T\leftrightarrow} \right)^k_{ij} = A^k_{ji} ;$$

$$\left( A^{T\searrow} \right)^{ij}_k = A^{kj}_i , \quad \left( A^{T\searrow} \right)^k_{ij} = A^j_{ik} ;$$

$$\left( A^{T\nearrow} \right)^{ij}_k = A^{ik}_j , \quad \left( A^{T\nearrow} \right)^k_{ij} = A^i_{kj} .$$

We shall use only a subset of these transpositions in the rest of this paper, but we have listed many options in order to provide a better perspective. In some cases, when we think it is clear what kind of transposition is appropriate, we shall use only a superscript $T$, in order to avoid too much clutter in the equations.

We will need some notation and terminology for vectors and co-vectors.
We will denote vectors (column vectors) by putting subscripts at the components of a vector, and superscripts on the components of a co-vector (row-vector).
The tensor product of two vectors $u = \left[ u_i \right]$ and $v = \left[ v_j \right]$, of possibly different dimensions, is defined as

$$u \otimes v = \left[ u_i v_j \right]$$

and is, therefore, a two-indexed array.

With a tridimensional matrix $A = \left[ A_i^{jk} \right]$, the bilinear operator $A(u \otimes v) = \left[ A_i^{jk} u_j v_k \right]$. The transposition of a vector $u = \left[ u_i \right]$ is defined as $u^T := \left[ u_i \, \delta^{ij} \right]$, where $\left[ \delta^{ij} \right]$ is an array of Kronecker's deltas, $\delta^{ii} = 1$, and, for $i \neq j$, $\delta^{ij} = 0$. Thus $\left( u^T \right)^i = u_i$. Likewise for co-vectors, if $w = \left[ w^j \right]$ is a co-vector, then its transpose is $w^T = \left[ w^j \, \delta_{ji} \right]$, in other words $\left( w^T \right)_i = w^j$.

A tridimensional matrix $A = \left[ A_i^{jk} \right]$ can act in 2 ways on a vector $w = \left[ w_i \right]$:

$$\overset{(1)}{A} w = \left[ \left( \overset{(1)}{A} w \right)_i^k \right] := \left[ A_i^{jk} w_j \right], \ \overset{(2)}{A} w = \left[ \left( \overset{(2)}{A} w \right)_i^j \right] := \left[ A_i^{jk} w_k \right].$$

The tridimensional matrix $A = \left[ A_i^{jk} \right]$ also defines a trilinear form. If $u = \left[ u_j \right]$, $v = \left[ v_k \right]$, $w = \left[ w^j \right]$, then the associated trilinear form is defined as

$wA(u \otimes v)$,

which equals $w^j A_i^{jk} u_j u_k$.

There are also three bilinear operations: $A(u \otimes v), \ w \overset{(1)}{A} u, \ w \overset{(2)}{A} v$.

We remark that, with the definitions we have formulated, we have certain identities about bilinear operations, for example, for a vector $u$ and a co-vector $w$, we have

$$w \overset{(1)}{A} u = u^T (A^{T+}) (\overset{(2)}{w^T}).$$

We include a few definitions and notation in differential calculus.

When $f$ is a differentiable function from an open set $G$ of $R^d$ into $R$, its gradient is a co-vector

$$\nabla_x f := \left[ \frac{\partial f}{\partial x_i} \right]_{1 \leq i \leq d} \text{ and its Hessian is a matrix (with both indices written as superscripts)}$$

$$\nabla_{xx} f := \left[ \frac{\partial^2 f}{\partial x_i \, \partial x_j} \right]_{\substack{1 \leq i \leq d \\ 1 \leq j \leq d}}.$$

Here, we are following one of the standard notational conventions of tensor analysis, viz. that a sub- (super-) script in a symbol that appears in a denominator is counted as a super- (sub-) script for the entire expression.

When the co-domain of $f$ is $R^n$, then its Jacobian matrix and its Hessian (tridimensional) matrix are

$$\nabla_x f = \left[\left(\nabla_x f\right)_i^j\right] := \left[\frac{\partial f_i}{\partial x_j}\right], \quad \nabla_{xx} f = \left[\left(\nabla_{xx} f\right)_i^{jk}\right] := \left[\frac{\partial^2 f_i}{\partial x_j \partial x_k}\right].$$

### 3. Controlled double Fredholm integral equations: first-order necessary conditions.

We consider a double Fredholm integral equation

$$\varphi(x) = f_0(x) + \int_G f_1(x, y, \varphi(y), u(y)) \, dy + \tfrac{1}{2} \int_G \int_G f_2(x, y, z, \varphi(y), \varphi(z), u(y), u(z)) \, dy \, dz$$

$$(3.1)$$

Here, $G$ is a bounded open set in real $d$ – dimensional Euclidean space, the state $\varphi$ takes values in real $n$ – dimensional space, and the control $u$ takes values in real $m$ – dimensional space. We assume that all functions are sufficiently many times continuously differentiable so that all indicated derivatives, appearing in our derivations below, exist and are continuous. We postulate existence and uniqueness of solution of (3.1).

The objective of an optimal control problem is the minimization of a cost functional

$$J := \int_G F_1(x, \varphi(x), u(x)) \, dx + \tfrac{1}{2} \int_G \int_G F_2(x, z, \varphi(x), \varphi(z), u(x), u(z)) \, dz \, dx \qquad (3.2)$$

The optimal control problem is non-Pontryaginian, the Hamiltonian is a functional of the state and the control, in addition to being a functional of the co-state. For single integral equations (degree of multiplicity equals 1), the Hamiltonian is a functional of the co-state but a pointwise function of the state and the control, and it is possible to prove an extremum principle of Pontryagin's type. For integral equations of degree of multiplicity > 1, the Hamiltonian is a functional also of the state and the control, and a pointwise extremum principle is not feasible.

We note that examples of various types of non-Pontryaginian control problems have been previously reported in several publications, for example in [W, AR]. A broad characterization of Pontryaginian optimal control problems is contained in [AD].

The optimal control theory of second-order integral equations is a natural extension of calculus of variations for multiple integrals, cf., e.g., [K, P].

We proceed to the derivation of the Hamiltonian equations, and, once we have the form of these equations, we prove rigorous results.

We assume that the functions $f_2$ and $F_2$ are symmetric, in the sense that

$$f_2(x, y, z, \varphi_1, \varphi_2, u_1, u_2) = f_2(x, z, y, \varphi_2, \varphi_1, u_2, u_1),$$
$$F_2(x, z, \varphi_1, \varphi_2, u_1, u_2) = F_2(z, x, \varphi_2, \varphi_1, u_2, u_1). \qquad (3.3)$$

This assumption entails no loss of generality, because, for arbitrary functions $f_2$ and $F_2$, we can replace $f_2$ and $F_2$ by their symmetrizations $\tilde{f}_2$ and $\tilde{F}_2$, respectively, where

$$\tilde{f}_2(x,y,z,\varphi_1,\varphi_2,u_1,u_2) := \tfrac{1}{2}[f_2(x,y,z,\varphi_1,\varphi_2,u_1,u_2) + f_2(x,z,y,\varphi_2,\varphi_1,u_2,u_1)],$$
$$\tilde{F}_2(x,z,\varphi_1,\varphi_2,u_1,u_2) := \tfrac{1}{2}[F_2(x,z,\varphi_1,\varphi_2,u_1,u_2) + F_2(z,x,\varphi_2,\varphi_1,u_2,u_1)]$$

(3.4)

without affecting the values of the corresponding integrals.

The co-state $\psi$ is a function having $G$ as its domain, and it takes values in the space of $d-$ dimensional real co-vectors (row vectors).
We define the penalty functional by

$$P := \int_G \psi(x) \Bigg[ \varphi_0(x) - \varphi(x) + \int_G f_1(x,y,\varphi(y),u(y))\,dy +$$
$$+ \tfrac{1}{2}\int_G \int_G f_2(x,y,z,\varphi(y),\varphi(z),u(y),u(z))\,dy\,dz \Bigg] dx$$

(3.5)

We denote by $\delta$ the operation of variations (in the sense of the Calculus of Variations), and by $\tilde{\delta}$ variations with respect to the state $\varphi$ only.
Towards the discovery of the appropriate Hamiltonian and the concomitant Hamiltonian equations, we set

$$\tilde{\delta}(P+J) = 0\,.$$

(3.6)

Taking into account the symmetry of the functions $f_2$ and $F_2$, and with some changes in the names of variables in various integrals, we find

$$\tilde{\delta}(P+J) = \int_G (-\psi(x))\,\delta\varphi(x)\,dx + \int_G \int_G \psi(y)\nabla_\varphi f_1(y,x,\varphi(x),u(x))\,dy\,\delta\varphi(x)\,dx +$$
$$+ \int_G \int_G \int_G \psi(y)\nabla_{\varphi(x)} f_2(y,x,z,\varphi(x),\varphi(z),u(x),u(z))\,dz\,dy\,\delta\varphi(x)\,dx +$$
$$+ \int_G \nabla_\varphi F_1(x,\varphi(x),u(x))\,\delta\varphi(x)\,dx + \int_G \int_G \nabla_{\varphi(x)} F_2(x,z,\varphi(x),\varphi(z),u(x),u(z))\,dz\,\delta\varphi(x)\,dx$$

(3.7)

Thus we are led to define the Hamiltonian

$$H(x,\varphi,\varphi(\cdot),u,u(\cdot),\psi(\cdot)) := F_1(x,\varphi,u) + \int_G F_2(x,z,\varphi,\varphi(z),u,u(z))\,dz +$$
$$+ \int_G \psi(y) f_1(y,x,\varphi,u)\,dy + \int_G \int_G \psi(y) f_2(y,x,z,\varphi,\varphi(z),u,u(z))\,dz\,dy$$

(3.8)

Then the condition

$\tilde{\delta}(P + J) = 0$

gives the Hamiltonian equation for the co-state:

$$\psi(x) = \nabla_\varphi H(x, \varphi, \varphi(\cdot), u, u(\cdot), \psi(\cdot)) \,. \qquad (3.9)$$

We shall prove:

<u>Theorem 3.1.</u> The total variation of the performance functional $J$ is given by

$$\delta J = \int_G (\nabla_{u(x)} H(x, \varphi(x), \varphi(\cdot), u(x), u(\cdot), \psi(\cdot))) \,\delta u(x) \, dx \,. \qquad (3.10)$$

<u>Proof</u>: From the definition of $J$, taking also into account the symmetry of $f_2$ and $F_2$, we obtain

$$\delta J = \int_G \{\nabla_\varphi F_1(x, \varphi(x), u(x)) \,\delta\varphi(x) + \nabla_u F_1(x, \varphi(x), u(x)) \,\delta u(x)\} \, dx +$$
$$+ \int_G \int_G \{\nabla_{\varphi(x)} F_2(x, z, \varphi(x), \varphi(z), u(x), u(z)) \,\delta\varphi(x) + \qquad (3.11)$$
$$+ \nabla_{u(x)} F_2(x, z, \varphi(x), \varphi(z), u(x), u(z)) \,\delta u(x)\} \, dz \, dx$$

From the equation of state dynamics, we have

$$\delta\varphi(x) = \int_G \{\nabla_\varphi f_1(x, y, \varphi(y), u(y)) \,\delta\varphi(y) + \nabla_u f_1(x, y, \varphi(y), u(y)) \,\delta u(y)\} \, dy +$$
$$+ \int_G \int_G \{\nabla_{\varphi(y)} f_2(x, y, z, \varphi(y), \varphi(z), u(y), u(z)) \,\delta\varphi(y) + \qquad (3.12)$$
$$+ \nabla_{u(y)} f_2(x, y, z, \varphi(y), \varphi(z), u(y), u(z)) \,\delta u(y)\} \, dz \, dy$$

Also, we write down, in explicit form, the Hamiltonian equation for the co-state:

$$\psi(x) = \int_G \psi(y) \nabla_\varphi f_1(y, x, \varphi(x), u(x)) \, dy + \int_G \int_G \psi(y) \nabla_{\varphi(x)} f_2(y, x, z, \varphi(x), \varphi(z), u(x), u(z)) \, dz \, dy +$$
$$+ \nabla_\varphi F_1(x, \varphi(x), u(x)) + \int_G \nabla_{\varphi(x)} F_2(x, z, \varphi(x), \varphi(z), u(x), u(z)) \, dz$$
$$\qquad (3.13)$$

We pre-multiply (3.12) by $\psi(x)$ and integrate over $G$ ; of course, in an integral, we can change the names of some variables as long as that change does not affect the value of the integral. Thus

$$\int_G \psi(x)\,\delta\varphi(x)\,dx + \delta J = \int_G \int_G \psi(y)\,\{\nabla_\varphi f_1(y,x,\varphi(x),u(x))\,\delta\varphi(x) +$$

$$+\nabla_u f_1(y,x,\varphi(x),u(x))\,\delta u(x)\}\,dy\,dx +$$

$$+\int_G \int_G \int_G \psi(y)\,\{\nabla_{\varphi(x)} f_2(y,x,z,\varphi(x),\varphi(z),u(x),u(z))\delta\varphi(x) +$$

$$+\nabla_{u(x)} f_2(y,x,z,\varphi(x),\varphi(z),u(x),u(z))\,\delta u(x)\}\,dz\,dy\,dx +$$

$$+\int_G \{\nabla_\varphi F_1(x,\varphi(x),u(x))\,\delta\varphi(x) + \nabla_u F_1(x,\varphi(x),u(x))\,\delta u(x)\}\,dx +$$

$$+\int_G \int_G \{\nabla_{\varphi(x)} F_2(x,z,\varphi(x),\varphi(z),u(x),u(z))\,\delta\varphi(x) +$$

$$+\nabla_{u(x)} F_2(x,z,\varphi(x),\varphi(z),u(x),u(z))\,\delta u(x)\}\,dz\,dx \qquad (3.14)$$

In view of the Hamiltonian equation, the integral $\int_G \psi(x)\,\delta\varphi(x)\,dx$ equals the sum of integrals of terms having a factor $\delta\varphi(x)$, and, after proper cancellations on the two sides of the last equality, we are left with

$$\delta J = \int_G \int_G \psi(y)\nabla_u f_1(y,x,\varphi(x),u(x))\,\delta u(x)\,dy\,dx +$$

$$+\int_G \int_G \int_G \psi(y)\,\nabla_{u(x)} f_2(y,x,z,\varphi(x),\varphi(z),u(x),u(z))\,\delta u(x)\,dz\,dy\,dx +$$

$$+\int_G \nabla_u F_1(x,\varphi(x),u(x))\,\delta u(x)\,dx + \qquad (3.15)$$

$$+\int_G \int_G \nabla_{u(x)} F_2(x,z,\varphi(x),\varphi(z),u(x),u(z))\,\delta u(x)\,dz\,dx$$

which is the explicit form of what was to be proved. ///

4. <u>The second variation for controlled double Fredholm integral equations.</u>

As mentioned before, an extremum principle of Pontryagin's type is not possible for double integral equations. Thus we look for alternatives to obtain more information, additional to the total variation formula obtained in the previous section. One possibility is the calculation of the second variation. It will turn out that the second variation requires, in addition to the Hamiltonian $H$ described above, a second Hamiltonian which we will term second-order ancillary Hamiltonian . The reasons for this terminology will become apparent as we proceed with the derivation of the equations of second variation.

We want to point out that our organization of the calculations is not the only possible; we have chosen to organize the work in a way that preserves as much as possible of the classical (first-order, in the sense of integral multiplicity) theory of optimal control of integral equations, and sets apart, in the most succinct way, the new ingredients that arise in second-order integral equations.

We evaluate the second variation of the state. We have

$$
\delta^2 \varphi(x) = \int_G \{ \nabla_\varphi f_1(x, y, \varphi(y), u(y)) \delta^2 \varphi(y) + \nabla_u f_1(x, y, ...) \delta^2 u(y) \} \, dy +
$$

$$
+ \int_G \{ \nabla_{\varphi\varphi} f_1(x, y, ...) (\delta\varphi(y) \otimes \delta\varphi(y)) + 2 \nabla_{\varphi u} f_1(x, y, ...) (\delta\varphi(y) \otimes \delta u(y)) +
$$

$$
+ \nabla_{uu} f_1(x, y, ...) (\delta u(y) \otimes \delta u(y)) \} \, dy +
$$

$$
+ \int_G \int_G \{ \nabla_{\varphi(y)} f_2(x, y, z, \varphi(y), \varphi(z), u(y), u(z)) \delta^2 \varphi(y) +
$$

$$
+ \nabla_{u(y)} f_2(x, y, z, \varphi(y), \varphi(z), u(y), u(z)) \delta^2 u(y) \} \, dy \, dz +
$$

$$
+ \int_G \int_G \{ \nabla_{\varphi(y)\varphi(y)} f_2(x, y, z, \varphi(y), \varphi(z), u(y), u(z)) (\delta\varphi(y) \otimes \delta\varphi(y)) +
$$

$$
+ 2 \nabla_{\varphi(y)u(y)} f_2(x, y, z, ...) (\delta\varphi(y) \otimes \delta u(y)) + \nabla_{u(y)u(y)} f_2(x, y, z, ...) (\delta u(y) \otimes \delta u(y)) +
$$

$$
+ \nabla_{\varphi(y)\varphi(z)} f_2(x, y, z, ...) (\delta\varphi(y) \otimes \delta\varphi(z)) + 2 \nabla_{\varphi(y)u(z)} f_2(x, y, z, ...) (\delta\varphi(y) \otimes \delta u(z)) +
$$

$$
+ \nabla_{u(y)u(z)} f_2(x, y, z, ...) (\delta u(y) \otimes \delta u(z)) \} \, dy \, dz
$$

(4.1)

We also calculate the initial expression for the second variation of $J$:

$$\delta^2 J = \int_G \ \{\nabla_\varphi F_1(x,\varphi(x),u(x))\,\delta^2\varphi(x) + \nabla_u F_1(x,\varphi(x),u(x))\,\delta^2 u(x)\}\,dx +$$

$$+\int_G \ \{\nabla_{\varphi\varphi}F_1(x,...)(\delta\varphi(x)\otimes\delta\varphi(x)) + 2\nabla_{\varphi u}F_1(x,...)(\delta\varphi(x)\otimes\delta u(x)) +$$

$$+\nabla_{uu}F_1(x,...)(\delta u(x)\otimes\delta u(x))\}\,dx +$$

$$+\int_G\int_G \ \{\nabla_{\varphi(x)}F_2(x,z,\varphi(x),\varphi(z),u(x),u(z))\,\delta^2\varphi(x) + \nabla_{u(x)}F_2(x,z,...)\,\delta^2 u(x)\}\,dz\,dx +$$

$$+\int_G\int_G \ \{\nabla_{\varphi(x)\varphi(x)}F_2(x,z,...)(\delta\varphi(x)\otimes\delta\varphi(x)) + 2\nabla_{\varphi(x)u(x)}F_2(x,z,...)(\delta\varphi(x)\otimes\delta u(x)) +$$

$$+\nabla_{u(x)u(x)}F_2(x,z,...)(\delta u(x)\otimes\delta u(x)) +$$

$$+\nabla_{\varphi(x)\varphi(z)}F_2(x,z,...)(\delta\varphi(x)\otimes\delta\varphi(z)) + 2\nabla_{\varphi(x)u(z)}F_2(x,z,...)(\delta\varphi(x)\otimes\delta u(z)) +$$

$$+\nabla_{u(x)u(z)}F_2(x,z,...)(\delta u(x)\otimes\delta u(z))\}\,dz\,dx$$

(4.2)

At this point, we introduce the second-order ancillary Hamiltonian

$$\text{`}\,h_2(x,z,\varphi_1,\varphi_2,u_1,u_2,\psi(\cdot)) := F_2(x,z,\varphi_1,\varphi_2,u_1,u_2) + \int_G \psi(y)f_2(y,x,z,\varphi_1,\varphi_2,u_1,u_2)\,dy$$

(4.3)

We shall prove:

<u>Theorem 4.1.</u> The second variation of the cost functional is given by

$$\delta^2 J = \int_G \ \{\nabla_{u(x)}H(x,\varphi(x),\varphi(\cdot\ ),u,u(\cdot),\psi(\cdot))\}\,\delta^2 u(x)\,dx +$$

$$+\int_G \ \{\nabla_{\varphi(x)\varphi(x)}H(x,...)(\delta\varphi(x)\otimes\delta\varphi(x)) + 2\,\nabla_{\varphi(x)u(x)}H(x,...)(\delta\varphi(x)\otimes\delta u(x)) +$$

$$+\nabla_{u(x)u(x)}H(x,...)(\delta u(x)\otimes\delta u(x))\}\,dx +$$

$$+\int_G\int_G \ \{\nabla_{\varphi(x)\varphi(z)}h_2(x,z,\varphi(x),\varphi(z),u(x),u(z),\psi(\cdot))(\delta\varphi(x)\otimes\delta\varphi(z)) +$$

$$+2\nabla_{\varphi(x)u(z)}h_2(x,z,...)(\delta\varphi(x)\otimes\delta u(x)) + \nabla_{u(x)u(z)}h_2(x,z,...)(\delta u(x)\otimes\delta u(z))\}\,dx\,dz$$

(4.4)

<u>Proof</u> : For convenience and greater clarity, we give a succinct but complete and precise verbal description of the procedure of proof. We evaluate the expression $\delta^2 J + \int_G \psi(x)\delta^2\varphi(x)\,dx$, using the RHS of (4.1) for $\delta^2\varphi(x)$, and the RHS of (4.2) for $\delta^2 J$. The terms involving a factor $\delta^2\varphi$ cancel from the two sides of the resulting equation, because of the Hamiltonian equation (3.13). The equation remaining, after such cancellations, is precisely the wanted equation (4.4). ///

<u>Remark 4.1.</u> The derivation of the equation of second variation is valid for general second-order variations of the control, i.e. the control $u$ is perturbed to

$$u_\varepsilon(x) := u(x) + \varepsilon\,\delta u(x) + \frac{\varepsilon^2}{2}\,\delta^2 u(x) + O(\varepsilon^3)\,.$$

For linear variations of the control, i.e. for $u_\varepsilon(x) = u(x) + \varepsilon\,\delta u(x)$, the second variation of the control vanishes, and the term that has $\delta^2 u$ as a factor, in (4.4), is absent in the final formula for $\delta^2 J$. Another case in which $\delta^2 u$ disappears is the case of no constraints on the control, so that, for an optimal control, $\nabla_u H = 0$. ///

5. <u>Sufficient conditions for optimality for Fredholm double integral equations.</u>

We consider a controlled linear Fredholm integral equation

$$\Phi(x) = \int_G \ [A(x,y)\Phi(y) + B(x,y)U(y)]\,dy \qquad (5.1)$$

and a cost functional

$$J_a := \int_G \ [\Phi^T(x)P_1(x)\Phi(x) + 2\Phi^T(x)Q_1(x)U(x) + U^T(x)R_1(x)U(x)]\,dx +$$

$$+ \int_G \int_G [\Phi^T(x)P_2(x,z)\Phi(z) + 2\Phi^T(x)Q_2(x,z)U(z) + U^T(x)R_2(x,z)U(z)\}\,dz\,dx \qquad (5.2)$$

We seek conditions under which $J_a$, under (5.1), is positive definite as a functional of $U$. Then, by applying this result to (3.) and (4.), with $\delta\varphi$, $\delta u$ in lieu of $\Phi$, $U$, respectively, we will obtain sufficient conditions of the Legendre-Clebsch type for the general problem of optimal control of Fredholm double integral equations.

We assume the unique solvability of (5.1), with a resolvent kernel $S(x,y)$, so that the solution of (5.1) can be represented in the form

$$\Phi(x) = \int_G \ B(x,y)U(y)\,dy + \int_G \int_G \ S(x,z)B(z,y)U(y)\,dz\,dy \equiv \int_G \ C(x,y)U(y)\,dy\,;$$

$$C(x,y) := B(x,y) + \int_G \ S(x,z)R(z,y)\,dz \qquad (5.3)$$

Merely for the purpose of being able to refer to these formulas, we write down the expression for the kernel $K$; we shall use the operator $Sym$ to denote the symmetrization of a square matrix valued function of two variables, namely

$$Sym(L(x,y)) := \tfrac{1}{2}(L(x,y) + L^T(y,x))\,.$$

Then

$$K(x_1,x_2) = \int_G \ C^T(x,x_1)P_1(x)C(x,x_2)\,dx + 2\,Sym\Big(C^T(x_2,x_1)Q_1(x_2)\Big) +$$

$$+ \int_G \int_G \ C^T(x,x_1)P_2(x,z)C(z,x_2)\,dx\,dz + 2\,Sym\Big(\int_G \ C^T(x,x_1)Q_2(x,x_2)\,dx\Big) + R_2(x_1,x_2)$$

Then $J_a$ takes the form

$$J_a = \int_G U^T(x) R_1(x) U(x)\, dx + \int_G \int_G U^T(x_1, x_2) K(x_1, x_2) U(x_2)\, dx_1\, dx_2 \qquad (5.4)$$

We say that the quadratic form, over $L^2(G)$, defined by the RHS of (5.4), is positive definite if its value is positive for every nonzero $U$ in $L^2(G)$.

We denote by $|G|$ the $d$ – dimensional Lebesgue measure of the domain $G$.

We define a kernel $M$ by

$$M(x_1, x_2) = \begin{bmatrix} \dfrac{R_1(x_1)}{|G|} & K(x_1, x_2) \\[2mm] K(x_2, x_1) & \dfrac{R_1(x_2)}{|G|} \end{bmatrix} \qquad (5.5)$$

Then

$$\int_G U^T(x) R_1(x) U(x)\, dx + \int_G \int_G U^T(x_1) K(x_1, x_2) U(x_2)\, dx_1\, dx_2 =$$
$$= \tfrac{1}{2} \int_G \int_G [U^T(x_1)\ \ U^T(x_2)] M(x_1, x_2) \begin{bmatrix} U(x_1) \\ U(x_2) \end{bmatrix} dx_1\, dx_2 \qquad (5.6)$$

Consequently, if the matrix $M(x_1, x_2)$ is positive definite for all $(x_1, x_2) \in \overline{G} \times \overline{G}$, then the functional $L_a$ is positive definite on $L^2(G \mapsto R^d)$.

6.  The quasi-LQC problem for second-order Fredholm integral equations.

These are systems that fall in the general category of the control systems studied in this paper, and lead to an integral variant of a Riccati equation.

We consider a system of second order of integral multiplicity, of the form

$$\varphi(x) = \varphi_0(x) + \int_G \ [f_0(x, y, \varphi(y)) + f_1(x, y, \varphi(y)u(y)]\,dy +$$

$$+ \tfrac{1}{2}\int_G \int_G \ [g_0(x, y, z, \varphi(y), \varphi(z)) + g_1(x, y, z, \varphi(y), \varphi(z))u(y)]\,dz\,dy \tag{6.1}$$

with a cost functional, also of the second order of integral multiplicity, and quadratic in the control:

$$J_2 := \int_G \ [F_0(x, \varphi(x)) + F_1(x, \varphi(x))u(x) + \tfrac{1}{2}u^T(x)F_2(x, \varphi(x))u(x)]\,dx +$$

$$+ \tfrac{1}{2}\int_G \int_G \ [F_3(x, z, \varphi(x), \varphi(z)) + F_4(x, z, \varphi(x), \varphi(z))u(x) + \tfrac{1}{2}u^T(x)F_5(x, z, \varphi(x), \varphi(z))u(x) +$$

$$+ \tfrac{1}{2}u^T(x)F_6(x, z, \varphi(x), \varphi(z))u(z)]\,dz\,dx \tag{6.2}$$

The Hamiltonian is

$$H(x, \varphi, \varphi(\cdot), u, u(\cdot), \psi(\cdot)) := F_0(x, \varphi) + F_1(x, \varphi)u + \tfrac{1}{2}u^T F_2(x, \varphi)u +$$

$$+ \int_G \ [F_3(x, z, \varphi, \varphi(z)) + F_4(x, z, \varphi, \varphi(z))u + \tfrac{1}{2}u^T F_5(x, z, \varphi, \varphi(z))u +$$

$$+ \tfrac{1}{2}u^T F_6(x, z, \varphi, \varphi(z))u(z)]\,dz + \int_G \ \psi(y)[f_0(y, x, \varphi) + f_1(y, x, \varphi)u]\,dy +$$

$$+ \int_G \int_G \ \psi(y)[g_0(y, x, z, \varphi, \varphi(z)) + g_1(y, x, z, \varphi, \varphi(z))u]\,dz\,dy \tag{6.3}$$

We assume no constraints on the control. We denote by $u^*$ a control that satisfies the condition

$$\nabla_u H(x, ..., u^*, u^*(\cdot), ...) = 0 \tag{6.4}$$

which, after a number of manipulations, leads to an integral equation of the form

$$u^{*T}(x) = G_0(x, \varphi(x), \varphi(\cdot)) + \int_G \ \psi(y)G_1(y, x, \varphi(x), \varphi(\cdot))\,dy + \int_G \ u^{*T}(z)G_2(z, x, \varphi(x), \varphi(\cdot))\,dz \tag{6.5}$$

The calculation of the precise expressions for $G_i$ , $0 \le i \le 2$ , is tedious but straightforward, and, for the purposes of the present paper, it is omitted. Here, we want to demonstrate the <u>form</u> of the relevant equations.

It may be convenient to also write down the transposed version of (6.5):

$$u^*(x) = G_0^T(x, \varphi(x), \varphi(\cdot)) + \int_G \ G_1^T(y, x, \varphi(x), \varphi(\cdot))\psi^T(y)\,dy + \int_G \ G_2^T(x, z, \varphi(x), \varphi(\cdot))u(z)\,dz$$

$$(6.6)$$

This is a linear Fredholm integral equation of the second kind, with kernel $G_2^T$. Assuming unique solvability, and the existence of a resolvent kernel, the solution has the representation

$$u^*(x) = S_0(x, \varphi(\cdot)) + \int_G \ S_1(x, y, \varphi(\cdot))\psi^T(y)\,dy \tag{6.7}$$

where, with $R(x, y)$ being the resolvent kernel,

$$S_0(x, \varphi(\cdot)) = G_0^T(x, \varphi(x), \varphi(\cdot)) + \int_G \ R(x, y)G_0^T(y, \varphi(y), \varphi(\cdot))\,dy \ ;$$

$$S_1(x, y, \varphi(\cdot)) = G_1^T(y, x, \varphi(x), \varphi(\cdot)) + \int_G \ R(x, z)G_1^T(y, z, \varphi(z), \varphi(\cdot))\,dz$$

$$(6.8)$$

The Hamiltonian equation for the co-state, $\psi(x) = \nabla_{\varphi(x)}H(x,...)$ , becomes

$$\psi(x) = \left(\nabla_{\varphi(x)}F_0(x,...)\right) + \left(\nabla_{\varphi(x)}F_1(x,...)\right)u^*(x) + \tfrac{1}{2}u^{*T}(x)\left(\nabla_{\varphi(x)}F_2(x,...)\right)u(x) +$$

$$+ \int_G \ \left[ \ \left(\nabla_{\varphi(x)}F_3(x, z,...)\right) + \left(\nabla_{\varphi(x)}F_4(x, z,...)\right)u^*(x) + \tfrac{1}{2}u^{*T}(x)\left(\nabla_{\varphi(x)}F_5(x, z,...)\right)u^*(x) +$$

$$+ \tfrac{1}{2}u^{*T}(x)\left(\nabla_{\varphi(x)}F_6(x, z,...)\right)u(z) \ \ \right]dz\,dx +$$

$$+ \int_G \ \psi(y)\nabla_{\varphi(x)}\left[ \ \ f_0(y, x,...) + f_1(y, x,...)u^*(x) + \int_G \ [g_0(y, x, z,...) + g_1(y, x, z,...)u^*(x)]dz \ \ \right]dy$$

$$(6.9)$$

### 7. The optimal control of second-order Volterra equations.

In this section, we shall briefly explain the c
corresponding problems of optimal control of double Volterra equations. Because the techniques we shall employ are similar to the case of double Fredholm equations, we shall omit many of the details, and we will be more explicit about those points that partly deviate from the Fredholm case.

A controlled double Volterra equation is

$$y(t) = y_0(t) + \int_0^t f_1(t,s,y(s),u(s))\,ds + \tfrac{1}{2}\int_0^t \int_0^t f_2(t,s,\sigma,y(s).y(\sigma),u(s),u(\sigma))\,d\sigma\,ds \qquad (7.1)$$

with cost functional, which is to be minimized,

$$J := F_0(T,y(T)) + \int_0^T F_1(t,y(t),u(t))\,dt + \tfrac{1}{2}\int_0^T\int_0^T F_2(t,\sigma,y(t),y(\sigma),u(t),u(\sigma))\,d\sigma\,dt \qquad (7.2)$$

For convenience, we also write down the equation for the value of the state at the final time $T$:

$$y(T) = y_0(T) + \int_0^T f_1(T,t,y(t),u(t))\,dt + \tfrac{1}{2}\int_0^T\int_0^T f_2(T,t,\sigma,y(t),y(\sigma),u(t),u(\sigma))\,d\sigma\,dt \qquad (7.3)$$

We postulate the usual conditions of differentiability of all involved functions, and conditions of symmetry of matrix-valued functions. We shall designate with $Y$ the variable that represents $y(T)$. We set

$$\omega := \nabla_Y F_0(T,y(T)) \qquad (7.4)$$

The Hamiltonian is defined as

$$\begin{aligned}
&H(t,T,y(t),y(T),y(\cdot),u(t),u(\cdot),\psi(\cdot),\omega) := \\
&= \omega\left[ f_1(T,t,y(t),u(t)) + \int_0^T f_2(T,t,\sigma,y(t),y(\sigma),u(t),u(\sigma))\,d\sigma \right] + \\
&+ F_0(T,y(T)) + F_1(t,y(t),u(t)) + \int_0^T F_2(t,\sigma,y(t),y(\sigma),u(t),u(\sigma))\,d\sigma + \\
&+ \int_t^T \psi(s) f_1(s,t,y(t),u(t))\,ds + \int_t^T\int_0^s \psi(s) f_2(s,t,\sigma,y(t),y(\sigma),u(t),u(\sigma))\,d\sigma\,ds
\end{aligned} \qquad (7.5)$$

The Hamiltonian equations are

$$\begin{aligned}
&\omega = \nabla_Y H(t,T,\ldots); \\
&\psi(t) = \nabla_{y(t)} H(t,T,\ldots)
\end{aligned} \qquad (7.6)$$

The formula of total variation is

$$\delta J = \int_0^T \left( \nabla_{u(t)} H(t, T, \ldots) \right) \delta u(t) \, dt \tag{7.7}$$

For the second variation of $J$, we shall utilize also the <u>ancillary second-order Hamiltonian</u> (analogous to the corresponding expression for the Fredholm case)

$$h_2(t, \sigma, T, y(t), y(\sigma), u(t), u(\sigma), \psi(\cdot)) := \omega f_2(T, t, \sigma, y(t), y(\sigma), u(t), u(\sigma)) +$$

$$+ F_2(t, \sigma, y(t), y(\sigma), u(t), u(\sigma)) + \int_{\max(t,\sigma)}^T \psi(s) f_2(s, t, \sigma, y(t), y(\sigma), u(t), u(\sigma)) \, ds \tag{7.8}$$

Then the second variation of $J$ takes the form

$$\delta^2 J = \int_0^T [\nabla_{u(t)} H(t, T, \ldots)] \delta^2 u(t) \, dt + (\nabla_{YY} F_0(T, y(T)))(\delta y(T) \otimes \delta y(T)) +$$

$$+ \int_0^T [(\nabla_{y(t)y(t)} H(t, T, \ldots))(\delta y(t) \otimes \delta y(t)) + 2(\nabla_{y(t)u(t)} H(t, T, \ldots))(\delta y(t) \otimes \delta u(t)) +$$

$$+ (\nabla_{u(t)u(t)} H(t, T, \ldots))(\delta u(t) \otimes \delta u(t))] \, dt +$$

$$+ \int_0^T \int_0^T [(\nabla_{y(t)y(\sigma)} h_2(t, \sigma, T, \ldots))(\delta y(t) \otimes \delta y(\sigma)) + 2(\nabla_{y(t)u(\sigma)} h_2(t, \sigma, T, \ldots))(\delta y(t) \otimes \delta u(\sigma)) +$$

$$+ (\nabla_{u(t)u(\sigma)} h_2(t, \sigma, T, \ldots))(\delta u(t) \otimes \delta u(\sigma))] \, d\sigma \, dt \tag{7.9}$$

We assume no constraints on the control; consequently, the term involving $\delta^2 u$, in (7.9), vanishes. In that case, therefore, the second variation is a quadratic functional, of the second order of integral multiplicity, in $\delta y$, $\delta u$. The dynamics of $(\delta y, \delta u)$ is a linear Volterra system of the form

$$\delta y(t) = \int_0^t \{A_1(t, s) \delta y(s) + B_1(t, s) \delta u(s)\} \, ds \tag{7.10}$$

and the <u>accessory</u> cost functional, i.e. the RHS of the equation for the second variation of $J$, is of the form

$$J_a = \frac{1}{2} (\delta y(T))^T P_0(T)(\delta y(T)) +$$

$$+ \int_0^T \{\tfrac{1}{2}(\delta y(t))^T P_1(t)(\delta y(t)) + (\delta y(t))^T Q_1(t)(\delta u(t)) + \tfrac{1}{2}(\delta u(t))^T R_1(t)(\delta u(t))\} \, dt +$$

$$+ \int_0^T \int_0^T \{\tfrac{1}{2}(\delta y(t))^T P_2(t, \sigma)(\delta y(\sigma)) + (\delta y(t))^T Q_2(t, \sigma)(\delta u(\sigma)) + \tfrac{1}{2}(\delta u(t))^T R_2(t, \sigma)(\delta u(\sigma))\} \, d\sigma \, dt \tag{7.11}$$

Our task now is to reduce (7.11) to the form

$$J_a = \frac{1}{2}\int_0^T (\delta u(t))^T R_1(t)(\delta u(t))\, dt + \frac{1}{2}\int_0^T \int_0^T (\delta u(t))^T K_2(t,\sigma)(\delta u(\sigma))\, d\sigma\, dt \qquad (7.12)$$

and then invoke the results of [BS] to establish conditions under which the pair of kernels $(R_1, K_2)$ determines a positive definite quadratic form on $L^2(G \rightarrow R^n)$. That will lead to sufficient conditions, of the Legendre - Clebsch type, for a solution of the set of first-order necessary conditions to be a local minimizer of $J$.

Starting with (7.10), we utilize a <u>resolvent kernel</u> $S(t,s)$ associated with the direct kernel $A_1(t,s)$, i.e. $S(t,s)$ solves

$$S(t,s) = A_1(t,s) + \int_s^t A_1(t,\sigma)S(\sigma,s)\, d\sigma \qquad (7.13)$$

Then the solution of (7.10) is

$$\delta y(t) = \int_0^t S_1(t,s)\,\delta u(s)\, ds \,; \quad S_1(t,s) := B_1(t,s) + \int_s^t S(t,\sigma)B_1(\sigma,s)\, d\sigma \qquad (7.14)$$

Utilizing (7.14) into (7.9) gives, after some integrations,

$$\begin{aligned}
K_2(t,\sigma) = {} & S_1^T(T,t)(\nabla_{yy} F_0(T,y(T))S_1(T,\sigma) + \\
& + \int_{\max(t,\sigma)}^T S_1^T(s,t)P_1(s)S_1(s,\sigma)\, ds + S_1^T(\sigma,t)Q_1(\sigma) + Q_1^T(t)S_1(t,\sigma) + \\
& + \int_\sigma^T \int_t^T S_1^T(s_1,t)P_2(s_1,s_2)S_1(s_2,\sigma)\, ds_1\, ds_2 + \\
& + \int_{\max(t,\sigma)}^T [S_1^T(s,t)Q_2(s,\sigma) + Q_2^T(s,t)S_1(s,\sigma)]\, ds + R_2(t,\sigma)
\end{aligned} \qquad (7.15)$$

## 8. Bilinear Volterra controlled systems of first and second orders of integral multiplicity.

### 8.1. First-order bilinear Volterra controlled system.

We consider the system

$$y(t) = y_0(t) + \int_0^t \ [A(t,s)\,y(s) + B(t,s)\,u(s) + y^T(s)\,C(t,s)\,u(s)]ds \tag{8.1}$$

with the optimal control objective of minimizing a functional

$$J := \int_0^T \ \{\tfrac{1}{2}\,y^T(t)P(t)y(t) + y^T(t)Q(t)u(t) + \tfrac{1}{2}u^T(t)R(t)u(t)\}\,dt \tag{8.2}$$

$C(t,s)$ is a tridimensional matrix (in the sense of indicial dimensionality), and the expression $y^T(s)C(t,s)u(s)$ has the same meaning as $C(t,s)(y(s) \otimes u(s))$.

The Hamiltonian is

$$H \equiv H(t,y,u,\psi(\cdot)) = \tfrac{1}{2}\,y^T P(t)y + y^T Q(t)u + \tfrac{1}{2}u^T R(t)u +$$
$$+ \int_t^T \ \psi(s)[A(s,t)y + B(s,t)u + y^T C(s,t)u]ds \tag{8.3}$$

The first-order necessary conditions for optimality are

$$\psi(t) = \nabla_y H \ ; \ \ \nabla_u H = 0 \tag{8.4}$$

thus

$$\psi(t) = y^T(t)P(t) + u^T(t)Q^T(t) +$$
$$+ \int_t^T \ \psi(s)[A(s,t) + u^T(t)C^T(s,t)]ds \ ; \tag{8.5}$$
$$y^T(t)Q(t) + u^T(t)R(t) + \int_t^T \ \psi(s)[B(s,t) + y^T(t)C(s,t)]ds = 0$$

The second equation in (8.4) gives

$$u^T(t) = -y^T(t)Q(t)R^{-1}(t) - \int_t^T \ \psi(s)[B(s,t) + y^T(t)C(s,t)]R^{-1}(t)\,ds \tag{8.6}$$

and substitution of (8.6) into (8.5) gives

$$\psi(t) = y^T(t)[P(t) - Q(t)R^{-1}(t)Q^T(t)] +$$

$$+(-1)\int_t^T \psi(s)[A(s,t) + B(s,t)R^{-1}(t)Q^T(t) + y^T(t)Q(t)R^{-1}(t)C^T(s,t)]ds +$$ (8.7)

$$+(-1)\int_t^T \int_t^T \psi(s_1)\psi(s_2)y^T(t)C(s_2,t)C^T(s_1,t)\,ds_2\,ds_1$$

For the second variation we have the equation

$$\delta^2 J = \int_0^T \left[ \delta y^T(t)P(t)\delta y(t) + 2\,\delta y^T(t)\left(Q(t) + \int_t^T \psi(s)C(s,t)\,ds\right)\delta u(t) + \delta u^T(t)D(t)\delta u(t) \right]dt$$

(8.8)

Consequently, the positive definiteness of the matrix

$$\begin{bmatrix} P(t) & Q(t) + \int_t^T \psi(s)C(s,t)\,ds \\ \left(Q(t) + \int_t^T \psi(s)C(s,t)\,ds\right)^T & D(t) \end{bmatrix}$$

(8.9)

is a sufficient condition for a solution of the first-order necessary conditions to be a local minimizer of $J$.

The sufficient condition can be made sharper, but in a more roundabout way, as follows: From the state dynamics, we have

$$\delta y(t) = \int_0^t \left\{ [A(t,s) + u^T(s)C^T(t,s)]\delta y(s) + [B(t,s) + y^T(s)C(t,s)]\delta u(s) \right\}ds$$ (8.10)

where $y$ and $u$ are solutions of the first-order necessary conditions. Then () is a second-kind linear Volterra equation and its solution is expressed, via a resolvent kernel, in the form

$$\delta y(t) = \int_0^t \Lambda(t,s)\,\delta u(s)\,ds$$ (8.11)

where the kernel $\Lambda$ is a function of $s$ and $t$ and a functional in $y, u$. (The actual calculation of $\Lambda$ is beyond the scope of the present discussion.) Substitution of (8.11) into () gives a quadratic integral expression for $\delta^2 J$ as a functional of $\delta u$, of the form

$$\delta^2 J = \int_0^T \delta u^T(t)\Lambda_1(t)\delta u(t)\,dt + \int_0^T \int_0^T \delta u^T(t)\Lambda_2(t,\tau)\delta u(\tau)\,d\tau\,dt$$

and then the second order sufficiency condition amounts to the requirement that the pair of kernels $(\Lambda_1, \Lambda_2)$ jointly determine a positive definite quadratic integral form on $L^2(0,T;R^m)$, i.e. the expression on the right-hand side of () above should be positive for all $\delta u \in L^2(0,T;R^m)$. The kernels $\Lambda_1$ and $\Lambda_2$ are given by

$$\Lambda_1(t) = R(t);$$
$$\Lambda_2(t,\tau) = Sym[\Lambda^T(\tau,t)(Q(\tau) + \int_\tau^T \psi(s)C(s,\tau)\,ds] + \int_{\max(t,\tau)}^T \Lambda^T(s,t)P(s)\Lambda(s,\tau)\,ds$$

## 8.2. The second-order bilinear Volterra controlled system.

The system has the form

$$y(t) = y_0(t) + \int_0^t [A(t,s)y(s) + B(t,s)u(s) + C(t,s)(y(s)\otimes u(s))]\,ds +$$
$$+ \int_0^t \int_0^t D(t,s,\sigma)(y(s)\otimes u(\sigma))\,d\sigma\,ds$$

and the associated optimal control problem is the minimization of the functional

$$J := \int_0^T [\tfrac{1}{2}y^T(t)P_1(t)y(t) + y^T(t)Q_1(t)u(t) + \tfrac{1}{2}u^T(t)R_1(t)u(t)]\,dt +$$
$$+ \int_0^T \int_0^T [\tfrac{1}{2}y^T(t)P_2(t,\tau)y(\tau) + y^T(t)Q_2(t,\tau)u(\tau) + \tfrac{1}{2}u^T(t)R_2(t,\tau)u(\tau)]\,d\tau\,dt$$

The Hamiltonian is

$$H = \tfrac{1}{2}y^T(t)P_1(t)y(t) + y^T(t)Q_1(t)u(t) + \tfrac{1}{2}u^T(t)R_1(t)u(t) +$$
$$+ \int_0^T [y^T(t)P_2(t,\tau)y(\tau) + y^T(t)Q_2(t,\tau)u(\tau) + u^T(t)Q_2^T(\tau,t)y(\tau) + u^T(t)R_2(t,\tau)u(\tau)]\,d\tau +$$
$$+ \int_t^T \psi(s)[A(s,t)y(t) + B(s,t)u(t) + C(s,t)(y(t)\otimes u(t))]\,ds +$$
$$+ \int_t^T \int_0^s \psi(s)[D(s,t,\sigma)(y(t)\otimes u(\sigma)) + D^{T\leftrightarrow}(s,\sigma,t)(u(t)\otimes y(\sigma))]\,d\sigma\,ds$$

In the unconstrained case, an optimal control is found from

$$(\nabla_u H)\big|_{u=u^*} = 0 \text{ which gives}$$

$$y^T(t)Q_1(t) + u^{*T}(t)R_1(t) +$$

$$+ \int_0^T [y^T(\tau)Q_2(\tau,t) + u^{*T}(\tau)R_2(\tau,t)]d\tau + \int_t^T \psi(s)[B(s,t) + C(s,t)\overset{(1)}{y}(t)]ds +$$

$$+ \int_t^T \int_0^s \psi(s)[D(s,\sigma,t)\overset{(1)}{y}(\sigma)]d\sigma\,ds = 0$$

The Hamiltonian equation for the costate is

$$\psi(t) = y^T(t)P_1(t) = u^T(t)Q_1^T(t) + \int_0^T [y^T(\tau)P_2(\tau,t) + u^T(\tau)Q_2^T(t,\tau)]d\tau +$$

$$+ \int_t^T \psi(s)[A(s,t) + C(s,t)\overset{(2)}{u}(t)]ds + \int_t^T \int_0^s \psi(s)D(s,t,\sigma)\overset{(2)}{u}(\sigma)d\sigma\,ds$$

The second variation of $J$ is given by

$$\delta^2 J = \int_0^T [\delta y^T(t)P_1(t)\delta y(t) + 2\delta y^T(t)\left(Q_1(t) + \int_t^T \psi(s)C(s,t)ds\right)\delta u(t) +$$

$$+ \delta u^T(t)R_1(t)\delta u(t)]dt + \int_0^T \int_0^T [\delta y^T(t)P_2(t,\tau)\delta y(\tau) +$$

$$+ 2\delta y^T(t)\left(Q_2(t,\tau) + \int_{\max(t,\tau)}^T \psi(s)D^T(s,t,\tau)ds\right)\delta u(\tau) + \delta u^T(t)R_2(t,\tau)\delta u(\tau)]d\tau\,dt$$

A sufficient condition for optimality is that (a) the first order necessary conditions are satisfied, and (b) the pair of matrices

$$M_1(t) = \begin{bmatrix} P_1(t) & Q_1(t) + \int_t^T \psi(s)C(s,t)ds \\ \left(Q_1(t) + \int_t^T \psi(s)C(s,t)ds\right)^T & R_1(t) \end{bmatrix},$$

$$M_2(t,\tau) = \begin{bmatrix} P_2(t,\tau) & Q_2(t,\tau) + \int_{\max(t,\tau)}^T \psi(s)D^T(s,t,\tau)ds \\ \left(Q_2(t,\tau) + \int_{\max(t,\tau)}^T \psi(s)D^T(s,t,\tau)ds\right)^T & R_2(t,\tau) \end{bmatrix}$$

jointly generate a positive definite quadratic form on $L^2(0,T;R^{n+m})$, in the sense that, for every nonzero $w \in L^2(0,T;R^{n+m})$ we have

$$\int_0^T w^T(t)M_1(t)w(t)dt + \int_0^T \int_0^T w^T(t)M_2(t,\tau)w(\tau)d\tau\,dt > 0.$$

Again, the second-order sufficiency condition can be made sharper by evaluating the first variation of the state and substituting into ():

$$\delta y(t) = \int_0^t \left\{ [A(t,s) + C(t,s)\overset{(2)}{u}(s)]\delta y(s) + [B(t,s) + C(t,s)\overset{(1)}{y}(s)]\delta u(s) \right\} ds +$$

$$+ \int_0^t \int_0^t [D(t,s,\sigma)\overset{(2)}{u}(\sigma)\delta y(s) + D(t,\sigma,s)\overset{(1)}{y}(\sigma)\delta u(s)]\, ds$$

from which it follows that $\delta y$ has the representation

$$\delta y(t) = \int_0^t M(t,s)\,\delta u(s)\, ds$$

By substituting the above expression for $\delta y$ into the equation for $\delta^2 J$, we obtain an expression for $\delta^2 J$ as a quadratic integral form in $\delta u$,

$$\delta^2 J = \int_0^t \delta u^T(t) M_1(t)\delta u(t)\, dt + \int_0^T \int_0^T \delta u^T(t) M_2(t,\tau)\delta u(\tau)\, d\tau\, dt$$

where

$$M_1(t) = R_1(t)$$

$$M_2(t,\tau) = R_2(t,\tau) + \int_{\max(t,\tau)}^T M^T(s,t)P_1(s)M(s,\tau)\, ds + \int_t^T \int_\tau^T M^T(s,t)P_2(t,\tau)M(\sigma,\tau)\, d\sigma\, ds +$$

$$+ Sym\left(M^T(\tau,t)\left(Q_1(\tau) + \int_t^T \psi(s)C(s,\tau)\, ds\right)\right) +$$

$$+ Sym\left(\int_t^T M^T(s,t)\left(Q_2(t,\tau) + \int_{\max(s,t)}^T \psi(\sigma)D^T(\sigma,s,\tau)\, d\sigma\right)ds\right)$$